\documentclass[11pt]{article}
\usepackage{graphicx,epsfig}
\usepackage{amsmath,amssymb,euscript}
\usepackage{latexsym}
\usepackage{color}
\newtheorem{theorem}{\bf Theorem}[section]

\newtheorem{lemma}[theorem]{\bf Lemma}
\newtheorem{proposition}[theorem]{\bf Proposition}

\newtheorem{corollary}[theorem]{\bf Corollary}

\newcommand{\Proofend}{\hfill$\diamondsuit$}
\newcommand\n{{\mathbf{n}}}
\newcommand\m{{\mathbf{m}}}
\newcommand\x{{\mathbf{x}}}

\def\NN{{\mathbb N}}
\def\RR{{\mathbb R}}
\def\N{{\mathbf N}}

\def\w{{\mathbf{w}}}
\def\PP{{\mathbb P}}
\def\C{{\mathcal{C}}}

\def\sspan{\textrm{span}}

\DeclareMathOperator{\IDC}{{IDC}}
\def\IDS{\textrm{IDS}}
\def\GDC{\textrm{GDC}}
\def\GDS{\textrm{GDS}}
\def\diag{{\mathrm{diag}}}

\def\al{{\alpha}}
\def\be{{\beta}}

\def\f{{\mathbf f}}

\def\S{{\mathbb S}}
\def\V{{\mathbf{V}}}
\def\w{{\mathbf{w}}}

\def\F{{\mathbf{F}}}

\def\NN{{\mathbb N}}
\def\RR{{\mathbb R}}
\def\N{{\mathbf N}}

\def\f{{\mathbf f}}

\def\W{{\mathbf{W}}}
\def\Y{{\mathbf{Y}}}
\def\X{{\mathbf{X}}}
\def\PP{{\mathbb P}}
\begin{document}
%\begin{frontmatter}
\title{Some remarks on filtered polynomial interpolation\\  at Chebyshev nodes}
\author{Donatella Occorsio and Woula Themistoclakis { \thanks{ corresponding author}}}
\maketitle

\begin{abstract}
The present paper concerns filtered de la Vall\'ee Poussin (VP) interpolation at the Chebyshev nodes of the four kinds. This approximation model is interesting for applications because it combines the advantages of the classical Lagrange polynomial approximation (interpolation and polynomial preserving) with the ones of filtered approximation (uniform boundedness of the Lebesgue constants and reduction of the Gibbs phenomenon). Here we focus on some additional features that are useful in the applications of  filtered VP interpolation. In particular, we analyze the simultaneous approximation provided by the derivatives of the VP interpolation polynomials. Moreover, we state the uniform boundedness of VP approximation operators  in some Sobolev and H\"older--Zygmund spaces where several integro--differential models are uniquely and stably solvable.
\end{abstract}

\section{Introduction}
Lagrange polynomial interpolation at Chebyshev nodes is widely used in the applications but, dealing with uniform norms, it is well known it does not provide an optimal approximation due to the unboundedness of the Lebesgue constants (see e.g. \cite{mastromilobook, trefethen}).

A common way to improve the approximation consists in applying some filter functions (see e.g. \cite{FiMhaPr, Fi.Th1, DeMa-Erb,Sloan, ThBa}), but in this way we generally lose the interpolation property. Nevertheless, in \cite{woula_99} it was first proved that if we apply a de la Vall\'ee Poussin (briefly VP) filter to the Fourier--Chebyshev partial sums of a given function $f$ and we discretize the associated Fourier coefficients by means of a suitable Gauss--Chebyshev quadrature rule, then we get filtered approximation polynomials that interpolate $f$ at the nodes of the applied quadrature rule.

Here we are going to deal with a generalization of the VP filtered interpolating polynomials introduced in \cite{woula_99}, which depends on the choice of two (integer) degree--parameters: $n$, that is the number of the interpolation nodes, and $0<m<n$ that is responsible for the action ray of the VP filter \cite{woula_NA, filtered1D}. The bivariate extension, via tensor product, of such filtered interpolation has been recently analyzed in \cite{paperoAMC, lncs_nostro} in comparison with Lagrange interpolation at the same nodes \cite{occorsiorusso2011}, as well as at the same number of Padua points \cite{Bos_Via, Bos_Xu, Caliari2008, Caliari2005}.

The present paper mainly focusses on the univariate case. For any kind of Chebyshev zeros of order $n$, $V_n^mf$ denotes the VP filtered polynomial that interpolates $f$ at these zeros. It is a polynomial of degree at most $n+m-1$ and coincides with $f$ in the case that $f$ is a polynomial of degree at most $n-m$.

As $n$ and $m$ increase, the behavior of the approximation $f\approx V_n^mf$  has been studied for general Jacobi weights in several papers (see e.g. \cite{woula_NA,ThL1} and the references therein). It is known that such  an approximation can be comparable with the best polynomial approximation whenever $n\sim m$ {\color{black}(i.e., $m< n\le c_1 m$ with $c_1>0$ independent of $n$ and $m$). For the particular Chebyshev case and any degree--parameters $n\sim m$, necessary and sufficient conditions for getting a near--best approximation error w.r.t. weighted uniform norm have been recently stated in \cite{filtered1D} by simple bounds on the Jacobi exponents of the norm's weight.} The Chebyshev case has been also analyzed in \cite{Ca.Th-wave, FisTh} from the wavelet point of view.

In this paper we aim to give some additional remarks on such kind of VP filtered interpolation at Chebyshev nodes.

The combination of the good approximation $f\approx V_n^mf$ with the interpolation and polynomial preserving properties, may allow to construct numerical (polynomial projection) methods for solving functional equations. To this aim it is important that the polynomial quasi--projection map $V_n^m:f\rightarrow V_n^mf$ is uniformly bounded, w.r.t. $n$ and $m$, in the couple of spaces where the equation is uniquely solvable \cite{CaCrJuLu00,OcDeBoPrandtl, operatoreD,air}. With regard to this problem we dedicate a section of the paper (Section 3) where we study the mapping properties of the VP projections in Sobolev and, more generally, H\"older--Zygmund spaces of locally continuous functions with uniform norms. The stated results have been recently applied in the construction of an efficient numerical method for solving Prandtl--type equations \cite{prandtl_nostro}.

In many contexts it is also useful to approximate the  function $f$ together with its first $r\ge 1$ derivatives, by using the same data for all the approximations. Such simultaneous approximation problems have been widely studied by many authors (see e.g. \cite{kilgore,kilg_sza, gonska, draga2015,dbocco,mastromilobook} and the references therein) and turn out to be useful in many applications. Recently, in \cite{kiani}, simultaneous approximation has been proposed in  deriving a numerical solver of ordinary differential equations applied  in the field of satellite geodesy.
In \cite{cao} a study of the approximation error by feedforward neural networks (FNNs)  for the simultaneous approximation of functions is proposed, it is especially interesting since neural networks are widely used in many applications including  image processing.

A typical approach to the simultaneous approximation question consists in taking a sufficiently good basic  approximation of $f$, whose derivatives are easy to compute and well--approximate the derivatives of $f$.

For instance, following this idea, in \cite{fil_oc_th} generalized Bernstein polynomials have been used to get quadrature rules at equidistant nodes for the simultaneous approximation of the Hilbert transform and its derivative, the Hadamard transform. Concerning the same problem, also Lagrange interpolating polynomials have been used to get quadrature rules of Hilbert and Hadamard transforms at Jacobi zeros (see \cite{mastronardi, dbocco,etna}).

In the present paper (Section 5) we are going to consider the simultaneous approximation of $f$ and its derivatives obtained by using a VP filtered interpolating polynomial as basic approximation tool. We will provide error estimates in the case $f$ belongs to some Sobolev type spaces and in the case $f\in C^q([-1,1])$, which is especially advisable in applications where the function $f$ doesn't have endpoint algebraic singularities. The stated results lay the foundation for future applications of the simultaneous approximation based on VP filtered interpolation.

Another feature  we aim to develop in this paper, regards the numerical fast computation of the VP interpolating polynomials and its derivatives. We will consider the univariate as well as the bivariate case, providing in both the cases simple formulas based on Fourier discrete block transforms that can be computed by means of fast algorithms (see e.g. \cite{tasche}). Moreover, in the bivariate case we estimate the VP interpolation error w.r.t. the weighted uniform norm in terms of some univariate errors of best approximation.

Finally, we focus on the problem of the optimal choice of the action ray $m$ of the VP filter we are using. By several numerical experiments it is known that, once fixed the number $n$ of interpolation nodes, suitable choices of the degree--parameter $m<n$ may improve the approximation (see e.g. \cite{woula_NA, ThBa, lncs_nostro, paperoAMC}).  Also it has been noticed that the Lebesgue constants results to be uniformly bounded {\color{black} w.r.t. the degree-parameters whenever we choose integers such that $m=\theta n$ where $0<\theta<1$ is a fixed, rational {\it localization parameter} and the numerical evidence shows a decreasing behavior of the Lebesgue constants w.r.t. $\theta$.} Here (Section 6) we look for  the theoretical reason to the numerical outputs of the performed numerical tests, by deriving estimates of the Lebesgue constants that make explicit  the dependence on the degree--parameters $n>m$ and hence on the localization parameter $\theta=m/n$.

The outline of the paper is the following. In Section 2 we recall some basic facts on filterd VP interpolation and focus on the computational aspects. Section 3 concerns the mapping properties, Section 4 deals with the bivariate VP filtered interpolation and Section 5 with the simultaneous approximation. Finally, Section 6 focusses on the dependence on the degree--parameters $m<n$ of the VP approximation $f\approx V_n^m f$.

\section{Basics}
Let us consider the four Chebyshev weights
\[
w_1(x):=\frac 1{\sqrt{1-x^2}},\quad w_2(x):=\sqrt{1-x^2},\quad
w_3(x):=\sqrt{\frac{1+x}{1-x}},\quad w_4(x):=\sqrt{\frac{1-x}{1+x}},
\]
and let us adopt the notation without subscript, namely $w(x)$, to mean anyone of them. Moreover, denote by $\{p_n(w,x)\}_n$ the associated system of orthonormal polynomials having, for any $n\in\NN$,  the following trigonometric form
\begin{equation}\label{pol}
p_n(w,x)=\left\{\begin{array}{ll}
%\sqrt{\frac 2\pi}\ T_n(x):=
\sqrt{\frac 2\pi}\ \cos[nt] \quad \left(\frac 1{\sqrt{\pi}}\ \mbox{for}\ n=0\right)&
\mbox{if $w=w_1$,}\\[.15in] \displaystyle
%\sqrt{\frac 2\pi}U_n(x):=
\sqrt{\frac 2\pi}\ \frac{\sin[(n+1)t]}{\sin t} \quad & \mbox{if $w=w_2$,}\\[.15in]\displaystyle
%\frac 1{\sqrt{\pi}}V_n(x):=
\frac 1{\sqrt{\pi}}\ \frac{\cos[(2n+1)t/2]}{\cos[t/2]}  \quad & \mbox{if $w=w_3$,}\\[.15in]\displaystyle
% \frac 1{\sqrt{\pi}}W_n(x):=
\frac 1{\sqrt{\pi}}\ \frac{\sin[(2n+1)t/2]}{\sin[t/2]}\quad & \mbox{if $w=w_4$,}
\end{array}\right.
\end{equation}
where $t=\arccos x\in [0,\pi]$,
being understood the  continuous extension in the not defined cases occurring when $t\in\{0,\pi\}$.

For any pair of positive integers $m<n$, the VP approximation polynomial $V_n^mf(x)$ of a given function $f$  at any $|x|\le 1$, is defined by (see e.g. \cite{filtered1D,woula_NA})
\begin{eqnarray}\label{VP}
V_n^mf(x)&:=&\sum_{k=1}^n f(x_k^n)\Phi_{n,k}^m(x),  \\
\label{fi-1}
\Phi_{n,k}^m(x)&:=&\lambda_{k}^n\sum_{j=0}^{n+m-1}\mu_{n,j}^m p_j(w,x_k^n)p_{j}(w,x),
\qquad k=1,\ldots,n,
\end{eqnarray}
where, for $k=1,\ldots, n$,
\begin{equation}\label{zeros}
x_k^n:=x_k^n(w)=\cos t_k^n, \qquad\mbox{with}\qquad t_k^n:=t_{k}^n(w)=\left\{\begin{array}{ll}
%\displaystyle
\frac{(2k-1)\pi}{2n} & \mbox{if $w=w_1$,}\\ [.15in]
%\displaystyle
\frac{k\pi}{n+1} & \mbox{if $w=w_2$,}\\ [.15in]
%\displaystyle
\frac{(2k-1)\pi}{2n+1} & \mbox{if $w=w_3$,}\\ [.15in]
%\displaystyle
\frac{2k\pi}{2n+1} & \mbox{if $w=w_4$,}
\end{array}\right.
%\qquad k=1,\ldots,n,
\end{equation}
are the zeros of $p_n(w,x)$,
\begin{equation}\label{lambda}
\lambda_{k}^n:=\lambda_{k}^n(w)= \frac 1{\displaystyle\sum_{j=0}^{n}p_j^2(w,x_k^n)}
=\left\{ \begin{array}{cl}
%\displaystyle
\frac \pi n &\mbox{if}\quad  w=w_1,\\[.5cm]
%\displaystyle
\frac{\pi}{n+1}\sin^2t_k^n &\mbox{if}\quad w=w_2,\\ [.5cm]
%\displaystyle
\frac{4\pi}{2n+1}\cos^2\frac{t_k^n}2  &\mbox{if}\quad w=w_3,\\ [.5cm]
%\displaystyle
\frac{4\pi}{2n+1}\sin^2\frac{t_k^n}2  &\mbox{if}\quad w=w_4,
\end{array}\right.
\end{equation}
are the related Christoffel numbers, and $\mu_{n,j}^m$ are the following VP filtering coefficients
\[
\mu_{n,j}^m:=\left\{\begin{array}{ll}
1 & \mbox{if}\quad j=0,\ldots, n-m,\\ [.1in]
\displaystyle\frac{n+m-j}{2m} & \mbox{if}\quad
n-m< j< n+m.
\end{array}\right.
\]
Using the Darboux kernels $K_r(x,y):=\sum_{j=0}^rp_j(w,x)p_j(w,y)$, it can be checked that the sum (\ref{fi-1}) equals the following delayed mean of Darboux kernels
\[
\Phi_{n,k}^m(x)=\frac{\lambda_{k}^n}{2m}\sum_{r=n-m}^{n+m-1}K_r(x_k^n, x),\qquad k=1,\ldots,n.
\]
Finally, another interesting form of the so--called fundamental VP polynomials $\Phi_{n,k}^m$ can be achieved by means of the polynomials \cite{Ca.Th-wave}
\begin{equation}\label{q-basis}
q_{n,j}^m(w,x):=\left\{\begin{array}{lll}
p_j(w,x)& & \mbox{if}\quad 0\le j\le n-m,\\ [.1in]
\gamma_{n,j}^mp_j(w,x)-\gamma_{n,2n-j}^m p_{2n-j}(w,x)
& & \mbox{if}\quad n\hspace{-.05cm}-\hspace{-.05cm}m<j<n,
\end{array}\right.
\end{equation}
where we set
\begin{equation}\label{gamma-q}
 \gamma_{n,j}^m:=\frac{m+n-j}{2m}.
 \end{equation}
These polynomials satisfy the orthogonality relation
\[
\int_{-1}^1 q_{n,j}^m(w,x)q_{n,i}^m(w,x)w(x)dx=0,\qquad \forall i\ne j\in\{0,\ldots, n-1\} .
\]
The following orthogonal expansion holds true
\begin{equation}\label{fi-sum}
\Phi_{n,k}^m(x)=\lambda_k^n\sum_{j=0}^{n-1} p_j(w,x_k^n)q_{n,j}^m(w,x),
\end{equation}
which yields
\begin{equation}\label{VP-q}
V_n^mf(x)=\sum_{j=0}^{n-1} c_{n,j}(f)q_{n,j}^m(w,x), \qquad c_{n,j}(f):=\sum_{k=1}^n \lambda_k^nf(x_k^n)p_j(w,x_k^n).
\end{equation}
In conclusion, we remark that in both the formulas (\ref{VP}) and (\ref{VP-q}), we need only $n$ data (that are the function values $\{f(x_k^n)\}_{k=1,..,n}$ or the discrete Fourier coefficients $\{c_{n,j}(f)\}_{j=0,..,n-1}$) for computing $V_n^mf(x)$ that is a polynomial of degree at most $n+m-1$.
\subsection{Computational Aspects}
Now we focus the attention on the fast computation of the VP polynomial $V_n^mf$ at a given set of $N$ points $\{z_\ell\in [-1,1] : \ \ell=1,\ldots,N\}$. To this aim, we recall the definitions of some Fourier Discrete block transforms of
 a given data set  (see e.g. \cite{tasche},  \cite[(A.8),(A.14)]{martucci})
 \[
 \mathcal{A}=\{a_{j,\ell}:\ j=0,\ldots,n-1, \ \ell=1,\ldots,N\}
 \]
\begin{itemize}
\item The Inverse Discrete  Cosine transform of $\mathcal{A}$:
    \[
    \IDC(\mathcal{A}):= \left\{b_{s,\ell}:=\sum_{j=0}^{n-1}w_j a_{j,\ell} \cos\left[\frac{j(2s-1)\pi}{2n}\right]: \ s=1,\ldots,n,\ \ell=1,\ldots,N\right\},
    \]
    being $w_j=\sqrt{\frac 1 n}$ if  $j=0$ and $w_j=\sqrt{\frac 2 n}$ otherwise.
    %if $1\le j<n$.
\item The Inverse Discrete  Sine transform of $\mathcal{A}$:
    \[
    \IDS(\mathcal{A}):= \left\{b_{s,\ell}:=
    \frac{2}{n+1}\sum_{j=0}^{n-1}a_{j,\ell}\sin\left[\frac{(j+1)s\pi}{n+1}\right] : \ s=1,\ldots,n,\ \ell=1,\ldots,N\right\}.
    \]
\item The Generalized Discrete  Cosine transform of $\mathcal{A}$:
    \[
    \GDC(\mathcal{A}):= \left\{b_{s,\ell}:=
    2\sum_{j=0}^{n-1}a_{j,\ell} \cos\left[\frac{(2j+1)(2s-1)\pi}{2(2n+1)}\right] : \ s=1,\ldots,n,\ \ell=1,\ldots,N\right\}.
    \]
\item The Generalized Discrete  Sine transform of $\mathcal{A}$:
    \[
    \GDS(\mathcal{A}):= \left\{b_{s,\ell}:=
    2\sum_{j=0}^{n-1}a_{j,\ell}\sin\left[\frac{(2j+1)s\pi}{2(n+1)}\right] : \ s=1,\ldots,n,\ \ell=1,\ldots,N\right\}.
    \]
    \end{itemize}
Setting
\[
\V :=\left[V_n^mf(z_1),\dots,V_n^mf( z_N)\right],\qquad\quad
\f := \left[f(x_1^n),\dots,f( x_n^n)\right],
\]
the following proposition collects the formulas to compute the approximation vector $\V\in\RR^N$ from the data vector $\f\in\RR^n$.
\begin{proposition}\label{prop-fft}
Let $\Lambda\in \RR^{n\times n}$ and $\mathcal{Q}\in\RR^{n\times N}$ be defined by
\[
\Lambda:=\diag\left[\sqrt{\lambda_k^n(w)}\right]_{k=1,..,n},\quad\mbox{and}\quad
\mathcal{Q}=\left[q_{n,j}^m(w,z_\ell)\right]_{0\le j< n, \ 1\le \ell \le N}.
\]
If ${\Phi}\in \RR^{n\times N}$ is given by the following transforms
\begin{eqnarray}
\label{fi1-trig} \mbox{Case $w=w_1$} &:&
{\Phi}= \IDC(\mathcal{Q}),\\
\label{fi2-trig} \mbox{Case $w=w_2$} &:&
{\Phi}=\sqrt{\frac{n+1}{2}}\ \IDS(\mathcal{Q}),\\
\label{fi3-trig} \mbox{Case $w=w_3$} &:&
{\Phi}=\frac{1}{\sqrt{2n+1}}\ \GDC(\mathcal{Q}),\\
\label{fi4-trig} \mbox{Case $w=w_4$} &:&
{\Phi}=\frac{1}{\sqrt{2n+1}}\  \GDS(\mathcal{Q}),
\end{eqnarray}
then we have
\begin{equation}\label{VP-vec}
\V=\f \cdot\Lambda\cdot  {\Phi},
\end{equation}
where $\cdot$ denotes the usual row-column product.
\end{proposition}
{\it Proof of Proposition \ref{prop-fft}}
By taking into account (\ref{VP}) it is obvious that (\ref{VP-vec}) holds if the entries of the matrix ${\Phi}$ are given by
\[
{\Phi}_{s,\ell}=\Phi_{n,s}^m(z_\ell),\qquad s=0,\ldots, n-1,\qquad \ell=1,\ldots,N.
\]
On the other hand, by (\ref{fi-sum}), (\ref{pol}) and (\ref{lambda}),  we get the following  "mixed-trigonometric" expressions of the VP fundamental polynomials
\[
\Phi_{n,k}^m(x)=\left\{\begin{array}{ll}
%\displaystyle
\frac{\sqrt{2\pi}}{n}\left[\frac 1 {\sqrt{2\pi}}+\sum_{j=1}^{n-1}\cos\left(j\frac{(2k-1)\pi}{2n}\right) q_{n,j}^m(w,x)\right],  & \mbox{if $w=w_1$},\\ [.15in]
%\displaystyle
\frac{\sqrt{2\pi}}{n+1}\sin \left( \frac{k\pi}{n+1}\right)\sum_{j=0}^{n-1}\sin\left((j+1)\frac{k\pi}{n+1}\right) q_{n,j}^m(w,x),  & \mbox{if $w=w_2$},\\ [.15in]
%\displaystyle
\frac{4\sqrt{\pi}}{2n+1}\cos \left(\frac{(2k-1)\pi}{2(2n+1)}\right)\sum_{j=0}^{n-1}\cos\left((2j+1)\frac{(2k-1)\pi}{2(2n+1)}\right) q_{n,j}^m(w,x), & \mbox{if $w=w_3$},\\ [.15in]
%\displaystyle
\frac{4\sqrt{\pi}}{2n+1}\sin \left(\frac{(2k-1)\pi}{2(2n+1)}\right)\sum_{j=0}^{n-1}\sin\left((2j+1)\frac{(2k-1)\pi}{2(2n+1)}\right) q_{n,j}^m(w,x), & \mbox{if $w=w_4$},
\end{array}\right.
\]
that means (\ref{fi1-trig})--(\ref{fi4-trig}) hold. \Proofend

\section{Mapping properties}
In order to construct the VP interpolating polynomial $V_n^mf$, it is sufficient that the function $f$ is defined at the Chebyshev nodes of order $n$. Here, we are going to focus on the case that $f$ is locally continuous on $[-1,1]$ (i.e., $f$ continuous in any $[a,b]\subset ]-1,1[$). In this way, the VP approximation $V_n^mf$ is defined for any pair of integers $n>m>0$ and we do not exclude that $f$ can be unbounded at the endpoints $\pm 1$. More precisely we suppose the behavior of $f$ at the endpoints is governed by a Jacobi weight
\[
u(x)=v^{\gamma,\delta}(x):=(1-x)^\gamma (1+x)^\delta, \qquad \gamma,\delta\ge 0,
\]
which is such that
\[
\lim_{x\rightarrow +1}f(x)u(x)=0\qquad \mbox{if $\gamma>0$},
\qquad \mbox{and}
\qquad
\lim_{x\rightarrow -1}f(x)u(x)=0\qquad \mbox{if $\delta>0$}.
\]
We denote by $C^0_u$ the space of all such functions equipped with the norm
\[
\|f\|_{C^0_u}:=\|fu\|= \max_{x\in [-1,1]}|f(x)u(x)|.
\]
It is well-known that Weierstrass approximation theorem holds in this Banach space and, denoting by $\PP_n$ the set of all algebraic polynomials of degree at most $n$, if we consider the error of best approximation of $f\in C^0_u$ in $\PP_n$, namely
\[
E_n(f)_u:=\inf_{P\in\PP_n}\|(f-P)u\|,
\]
then we have
\[
f\in C^0_u \ \Longleftrightarrow \ \lim_{n\rightarrow \infty} E_n(f)_u=0.
\]
In the sequel, we suppose the Jacobi weight $u=v^{\gamma,\delta}$ satisfies the following conditions
\begin{equation}\label{hp-u}
\begin{array}{llllll}
\diamond \ \mbox{Case $w=w_1$:} &  0\le \gamma\le 1 & \mbox{and} & 0\le \delta\le 1 && \\
\diamond\ \mbox{Case $w=w_2$:} & 0< \gamma\le 3/2 & \mbox{and} & 0< \delta\le 3/2 & \mbox{and} & -1\le\gamma-\delta \le 1\\
\diamond\ \mbox{Case $w=w_3$:} & 0\le \gamma\le 1 & \mbox{and} & 0< \delta\le 3/2 & \mbox{and} & \gamma-\delta \le 1/2\\
\diamond\ \mbox{Case $w=w_4$:} & 0< \gamma\le 3/2 & \mbox{and} & 0\le \delta\le 1 & \mbox{and} & \gamma-\delta \ge -1/2.
\end{array}
\end{equation}
It has been recently proved \cite{filtered1D} that these bounds are necessary and sufficient in order that {\color{black} the VP polynomial $V_n^m f$ is a near--best approximation polynomial of any $f\in C^0_u$, for any pair of integers $n,m\in\NN$ satisfying $m<n\le c_1 m$, with $c_1>1$ independent of $n,m$. In fact, using the notation $n\sim m$ to indicate that $n,m\in \NN$ are related as above, we have the following}
\begin{theorem}\label{th-C0}
 Let $n,m\in\NN$ be such that {\color{black} $n\sim m$.} The conditions (\ref{hp-u}) are necessary and sufficient in order that the operator $V_n^m:f\rightarrow V_n^mf$, considered as a map from $C^0_u$ into itself, is uniformly bounded w.r.t. $n$. Moreover, this is equivalent to having
\begin{equation}\label{err-VP-C0}
E_{n+m-1}(f)_u\le \|(V_n^m f-f)u\|\le \C E_{n-m}(f)_u, \qquad \forall f\in C^0_u,
\end{equation}
with $\C>0$ independent of $f$ and $n$.
\end{theorem}
By virtue of this theorem, for suitable choices of $u$, we have
\[
\lim_{\begin{array}{c}
n\rightarrow\infty\\
{\color{black}n\sim m}
\end{array}}
\|(V_n^mf-f)u\|=0,\qquad \forall f\in C^0_u.
\]
By virtue of (\ref{err-VP-C0}), the rate of convergence  of such error is comparable with that of the error of best approximation of $f$ by polynomials, which, in turn, depends on the smoothness of $f$.

In particular, if we consider the Sobolev type spaces of order $r\in\NN$
\[W_r(u):=\left\{f\in C^0_u: \mbox{$f^{(r-1)}$ is locally absolutely continuous and } \|f^{(r)} \varphi^{r}u\|<\infty\right\},\]
equipped with the  norm
\[\|f\|_{W_r(u)}:=\|fu\| + \|f^{(r)}\varphi^{r}u\|,\]
it is well-known that (see e.g. \cite[Th. 2.1.1]{DT})
\begin{equation}\label{Favard-Sob}
E_n(f)_u\leq \frac{\C}{n^r} \|f\|_{W_r(u)}, \quad \forall f\in W_r(u), \quad \forall r\in\NN,\qquad \C\neq \C(n,f),
\end{equation}
where here and throughout the paper  $\C$ denotes a positive constant having different meaning in different formulas and we  write $\C \neq \mathcal{C}(n,f,\ldots)$  to say that $\C$ is  independent of  $n,f,\ldots$.

Besides the convergence result stated by Theorem \ref{th-C0}, one of the attractive features of the VP operator $V_n^m$ is that it is a projection map onto the polynomial space
\[
S_n^m:=\sspan\left\{ q_{n,j}^m(w,x): \ j=0,\ldots,n-1\right\}, \qquad \forall n>m,
\]
which is usually called VP space and satisfies the nesting property
\[
\PP_{n-m}\subset S_n^m\subset \PP_{n+m-1}, \qquad \forall n>m.
\]
The polynomial preserving property
\begin{equation}\label{inva}
V_n^m P=P,\qquad \forall P\in S_n^m,\qquad n>m,
\end{equation}
allows to construct projection methods to find the numerical solution of any functional equation uniquely solvable in $C^0_u$. Nevertheless, supposing that such equations have a unique and stable solution only in some subspaces $X$ of $C^0_u$, {\color{black} it is desirable that the map $V_n^m:X\rightarrow X$ is  uniformly bounded  w.r.t. $n\sim m$ in order to get optimal error estimates}.

In the case that $X$ is a Sobolev subspace $W_r(u)$, the following theorem is fundamental in order to prove the convergence of numerical methods based on VP projection operators.
\begin{theorem}\label{th-Sob}
Let the bounds (\ref{hp-u}) be satisfied and let $n,m\in\NN$ be such that
{\color{black} $n\sim m$}. For any $r\in\NN$, the map $V_n^m:W_r(u)\rightarrow W_r(u)$ is uniformly bounded w.r.t. $n$. Moreover, the following error estimates hold for all $f\in W_r(u)$ and any $s\in\NN$ with $s\le r$
\begin{eqnarray}\label{err-Sob1}
\|(V_n^m f-f)u\|&\le& \frac{\C}{n^r}\ \|f\|_{W_r(u)}, \qquad \C\ne \C(n,f),\\
\label{err-Sob2}
\|V_n^m f-f\|_{W_s(u)}&\le& \frac{\C}{n^{r-s}}\ \|f\|_{W_r(u)}, \qquad \C\ne \C(n,f).
\end{eqnarray}
\end{theorem}
{\it Proof of Theorem \ref{th-Sob}}

The estimate (\ref{err-Sob1}) follows from (\ref{err-VP-C0})  and (\ref{Favard-Sob}), by taking into account that $n-m=(1-\theta)n$.

In order to prove (\ref{err-Sob2}), we recall that the Favard type estimate (see e.g.\cite[p.171]{mastromilobook})
\begin{equation}\label{Favard}
E_n(f)_v\le \frac{\C}{n}E_{n-1}(f')_{v\varphi}, \qquad \C\ne\C(n,f),
\end{equation}
holds for any Jacobi weight $v$ with nonnegative exponents.
Moreover, for any polynomial $P\in\PP_\mu$ the following estimate holds (see e.g. \cite[(4.3.17)]{mastromilobook})
\begin{equation}\label{der}
\|(f-P)^{(s)}u\varphi^s\|\le \C \left[\mu ^s \|(f-P)u\|+ E_{\mu -s}(f^{(s)})_{u\varphi^s}\right], \qquad \C\ne \C(f,P,\mu).
\end{equation}
By taking in (\ref{der}) $P=V_n^mf\in \PP_{n+m-1}$, using (\ref{err-Sob1}), applying (\ref{Favard}) $(r-s)$ times, and recalling that $m=\theta n$, we get
\begin{eqnarray*}
\|(f-V_n^mf)^{(s)}u\varphi^s\|&\le& \C \left[(n+m-1)^s \|(f-V_n^mf)u\|+ E_{n+m-1 -s}(f^{(s)})_{u\varphi^s}\right]\\
&\le& \C \frac{(n+m-1)^s}{n^r}\ \|f\|_{W_r(u)}+\frac{\C}{(n+m-1-r)^{r-s}}E_{n+m-1-r}(f^{(r)})_{u\varphi^r}\\
&\le& \frac{\C}{n^{r-s}}\|f\|_{W_r(u)}+ \frac{\C}{n^{r-s}}\|f^{(r)} u\varphi^r\|\\
&\le& \frac{\C}{n^{r-s}}\|f\|_{W_r(u)}.
\end{eqnarray*}
It follows that (\ref{err-Sob2}) holds for all positive integers $s\le r$ since, by the previous arguments, we obtain
\[
\|f-V_n^mf\|_{W_s(u)}= \|(f-V_n^mf)u\|+ \|(f-V_n^mf)^{(s)}u\varphi^s\|\le \frac{\C}{n^{r-s}}\|f\|_{W_r(u)}.
\]
Finally, by applying (\ref{err-Sob2}) with $s=r$, we have
\[
\|V_n^mf\|_{W_r(u)}\le \|V_n^mf - f\|_{W_r(u)}+ \|f\|_{W_r(u)}\le \C \|f\|_{W_r(u)},
\qquad \C\ne \C(n,f),
\]
which yields the uniform boundedness of the map $V_n^m:W_r(u)\rightarrow W_r(u)$ w.r.t. $n$.
\Proofend
\vspace{.5cm}\newline
In order to generalize (\ref{Favard-Sob}) to any real $r>0$, we introduce the following spaces
\[
Z_r (u):=\{f\in C^0_u \ :\ \sup_{n>0}(n+1)^rE_n(f)_u<\infty\}, \qquad r>0,
\]
equipped with the following norm
\[
\|f\|_{Z_r(u)}:=\|f u\|+\sup_{n>0}(n+1)^rE_n(f)_u, \qquad r>0.
\]
Obviously, by this definition, we have
\begin{equation}\label{Favard-Zig}
E_n(f)_u\leq \frac{\|f\|_{Z_r(u)}}{n^r} ,\qquad \forall f\in Z_r(u),\quad r>0.
\end{equation}
The spaces $Z_r(u)$ are also known as H\"older--Zygmund spaces and constitute a particular case of the Besov-type spaces studied in \cite{DT_Besov}, where they have been equivalently defined in terms of Ditzian--Totik moduli of smoothness (\cite[Th. 2.1]{DT_Besov}).

The following theorem states the uniform boundedness of a VP projection operator in H\"older--Zygmund spaces and generalizes to these spaces all the estimates previously stated in Theorem \ref{th-Sob} for Sobolev spaces.

\begin{theorem}\label{th-Zig}
Let the bounds (\ref{hp-u}) be satisfied and let $n,m\in\NN$ be such that {\color{black} $n\sim m$.} For any $r>0$, the map $V_n^m:Z_r(u)\rightarrow Z_r(u)$ is uniformly bounded w.r.t. $n$. Moreover, the following error estimates hold for all $f\in Z_r(u)$ and for any $0<s\le r$
\begin{eqnarray}\label{err-Zig1}
\|(V_n^m f-f)u\|&\le& \frac{\C}{n^r}\ \|f\|_{Z_r(u)}, \qquad \C\ne \C(n,f),\\
\label{err-Zig2}
\|V_n^m f-f\|_{Z_s(u)}&\le& \frac{\C}{n^{r-s}}\ \|f\|_{Z_r(u)}, \qquad \C\ne \C(n,f).
\end{eqnarray}
\end{theorem}
{\it Proof of Theorem \ref{th-Zig}}.

The estimate (\ref{err-Zig1}) is an immediate consequence of (\ref{err-VP-C0}) and (\ref{Favard-Zig}), as well as the uniform boundedness of the map $V_n^m: Z_r(u)\rightarrow Z_r(u)$ easily follows from (\ref{err-Zig2}) with $r=s$. Hence the crucial point is the proof of (\ref{err-Zig2}). \newline Taking into account that, by (\ref{err-Zig1}), we get
\begin{eqnarray}\nonumber
\|V_n^m f-f\|_{Z_s(u)}&=& \|(f-V_n^mf)u\|+\sup_{k>0}(k+1)^s E_k(f-V_n^mf)_u\\
\label{eq1}
&\le&
\frac{\C}{n^r}\ \|f\|_{Z_r(u)}+\sup_{k>0}(k+1)^s E_k(f-V_n^mf)_u,
\end{eqnarray}
we have to investigate only the last term of the above formula. Recalling that $V_n^mf\in \PP_{n+m-1}$, we observe that
\[
E_k(f-V_n^mf)_u=\inf_{P\in\PP_k}\|(f-V_n^mf-P)u\|
 \left\{\begin{array}{ll}
= E_k(f) & \mbox{if $k\ge n+m-1$,}\\
\le \|(f-V_n^mf)u\| & \mbox{if $k< n+m-1$.}
\end{array}\right.
\]
Consequently, by (\ref{Favard-Zig}) we get
\begin{eqnarray*}\nonumber
\sup_{k\ge n+m-1}(k+1)^s E_k(f-V_n^mf)_u &=& \sup_{k\ge n+m-1}(k+1)^s E_k(f)_u\\
\nonumber&\le& \sup_{k\ge n+m-1}\frac{(k+1)^s}{k^r}\ \|f\|_{Z_r(u)}\\
\label{eq2}&=& \frac{\C}{n^{r-s}}\ \|f\|_{Z_r(u)},
\end{eqnarray*}
and by (\ref{err-Zig1}) we have
\begin{eqnarray}\nonumber
\sup_{k< n+m-1}(k+1)^s E_k(f-V_n^mf)_u &\le & \|(f-V_n^mf)u\|
\sup_{k< n+m-1}(k+1)^s\\
\nonumber &\le& \C n^s  \|(f-V_n^mf)u\|\\
\label{eq3}
&\le& \frac{\C}{n^{r-s}}\ \|f\|_{Z_r(u)}.
\end{eqnarray}
Thus, the statement follows from the above estimates (\ref{eq1})--(\ref{eq3}). \Proofend
\section{On the bivariate case}
In \cite{lncs_nostro} multivariate filtered polynomials, associated with any filter function and any Jacobi weight, have been studied. The special (interpolating) case based on VP filter functions and the Chebyshev weight of first kind has been also investigated in \cite{paperoAMC} where the pros and cons have been analyzed in comparison with Lagrange interpolation at the same Chebyshev tensor--product grid as well as at the same number of Padua points \cite{Caliari2005}. Here, for completeness, we are going to consider  bivariate VP filtered polynomials based on anyone of the four kinds of the Chebyshev weights.

These polynomials can be easily deduced, via tensor product, from the univariate case. More precisely, for any degree--parameters
\[
\n:=(n_1, n_2)\in\NN^2,\quad \m:=(m_1, m_2)\in\NN^2,\qquad\mbox{with $n_i>m_i$, $i=1,2$}
\]
and for any couple of Chebyshev weights
\[
v_i(x):=(1-x)^{\alpha_i}(1+x)^{\beta_i},\qquad |\alpha_i|=|\beta_i|=\frac 12, \qquad i=1,2,
\]
we define the following bivariate fundamental VP polynomials
\begin{equation}\label{VPfund-2D}
\Phi_{\n,k,h}^{\m}(x,y):=\Phi_{n_1,k}^{m_1}(x)\Phi_{n_2,h}^{m_2}(y), \qquad k=1,\ldots, n_1, \quad h=1,\ldots, n_2,
\end{equation}
being $\Phi_{n_i,k}^{m_i}$ the unicvariate fundamental VP polynomials associated with the Chebyshev weight $v_i$, for $i=1,2$.

Hence, for any locally continuous function $f$, given the $n_1\times n_2$ values attained by $f$ at the point set
\begin{equation}\label{grid}
\X^\n:=\left\{\x_{k,h}^{\n}:= \left(x_k^{n_1}(v_1),\ x_h^{n_2}(v_2)\right):\quad
k=1,\ldots, n_1,\ h=1,\ldots, n_2
\right\},
\end{equation}
the bivariate VP filtered polynomial of $f$ is given by
\begin{equation}\label{VP-biv}
\V_\n^\m f(x,y)=
\sum_{k=1}^{n_1}\sum_{h=1}^{n_2}f(\x_{k,h}^{\n})\Phi_{\n,k,h}^{\m}(x,y),\qquad (x,y)\in [-1,1]^2.
\end{equation}
Similarly to the one dimensional case, $\V_\n^\m f$ interpolates $f$ at the Chebyshev grid $\X^\n$, i.e. we have the identities
\begin{equation}\label{interp-2D}
\V_\n^\m f(\x_{k,h}^{\n})=f(\x_{k,h}^{\n}), \qquad k=1,\ldots n_1,\quad h=1,\ldots, n_2.
\end{equation}
Moreover, set
\[
\S_\n^\m:= S_{n_1}^{m_1}(v_1)\otimes S_{n_2}^{m_2}(v_2)=\sspan\left\{{\Phi}_{\n,k,h}^\m(x,y): \ k=1,\ldots,n_1,\ h=1,\ldots, n_2\right\},
\]
we have the polynomial preserving property
\[
\V_\n^\m f=f,\qquad \forall f\in\S_\n^\m.
\]
As regards the fast computation of such bivariate VP polynomial at an arbitrary set of $N_1\times N_2$ points
\[
\Y^\N:=\left\{(a_\ell^{N_1}, \ b_k^{N_2}): \ \ell=1,\ldots, N_1,\quad k=1,\ldots,N_2\right\}
\]
setting
\begin{eqnarray*}
\W&:=&\left[\V_\n^\m f(a_i^{N_1}, \ b_j^{N_2})\right]_{1\le i\le N_1,\  1\le j\le N_2}
\\ [.1in]
\F &:=&\left[f(x_i^{n_1},x_j^{n_2})\right]_{1\le i\le n_1,\ 1\le j\le n_2},
\end{eqnarray*}
 and
$${\Phi}_{1}:=\left[\Phi_{n_1,k}^{m_1}(a_\ell^{N_1})\right]_{1\le k\le n_1, \ 1\le \ell \le N_1} ,$$
$${\Phi}_{2}:=\left[\Phi_{n_2,k}^{m_2}(b_\ell^{N_2})\right]_{1\le k\le n_2, \ 1\le \ell \le N_2} ,$$
we have the following matrix decomposition for the values of the VP polynomial $\V_\n^\m f$ at the point set $\Y^N$
$$\W={{\Phi}_{1}}^T \Lambda \F\Lambda {\Phi}_{2}.$$
By the arguments used in the one dimensional case, the matrices ${\Phi}_{i}$ can be achieved by means of fast algorithms that compute the following discrete sine or cosine transforms (cf. Subsection 2.1)
\begin{eqnarray}
\label{fi1} w=w_1 &\Longrightarrow&
{\Phi}_{i}=\IDC(\mathcal{Q}_i),\qquad i=1,2,\\
\label{fi2} w=w_2 &\Longrightarrow&
{\Phi}_{i}=\sqrt{\frac{n_i+1}{2}}\ \IDS(\mathcal{Q}_i),\qquad i=1,2,\\
\label{fi3} w=w_3 &\Longrightarrow&
{\Phi}_{i}=\frac{1}{\sqrt{2n_i+1}}\ \GDC(\mathcal{Q}_i),\qquad i=1,2,\\
\label{fi4} w=w_4 &\Longrightarrow&
{\Phi}_{i}=\frac{1}{\sqrt{2n_i+1}}\  \GDS(\mathcal{Q}_i),\qquad i=1,2,
\end{eqnarray}
of the matrices
\[
\mathcal{Q}_{1}:=\left[ q_{n_1,i}^{m_1}(v_1,a_j^{N_1})\right]_{1\le i\le n_1, \ 1\le j \le N_1}, \]
\[
\mathcal{Q}_{2}:=\left[ q_{n_2,i}^{m_2}(v_2, b_j^{N_2})\right]_{1\le i\le n_2, \ 1\le j \le N_2}.
\]
Concerning the approximation $f(x,y)\approx \V_\n^\m f(x,y)$, near--best error estimates have been stated in \cite{lncs_nostro} for general (even not interpolating) VP filtered polynomials related to any bivariate Jacobi weight $\w(x,y)=v_1(x)v_2(x)$. Moreover, the necessary and sufficient results stated for the particular 1st kind Chebyshev weight (cf.  \cite[Th. 3.1]{paperoAMC}) can be easily extended to any kind of Chebyshev weight by taking into account the respective univariate results stated by Theorem \ref{th-C0}.

Here, we are going to investigate the dependence of the bivariate error estimates on the smoothness properties that the function $f$ has w.r.t. a single variable.
More precisely, for any bivariate function $f$ and for any $|x|<1$, we consider the univariate functions obtained by keeping fixed at $x$ the first or the second variable, namely,  we define the following univariate functions of $t$
\[
\overline{f}_x(t):=f(x,t), \qquad \overline{\overline{f}}_x(t):=f(t,x), \qquad \forall |t|<1,
\]
having used the one/two bars notation to indicate the position (first/second) of the fixed variable $x$.

In order to include the case that $f$ may be unbounded on the boundary of the square $[-1,1]^2$, we consider the weights
\begin{equation}\label{ui}
u_i(x):=(1-x)^{\gamma_i}(1+x)^{\delta_i},\qquad \gamma_i,\delta_i\ge 0, \qquad i=1,2,
\end{equation}
and suppose that
\begin{equation}\label{hpf}
\overline{f}_x\in C^0_{u_2} \quad\mbox{and}\quad \overline{\overline{f}}_x\in C^0_{u_1}, \qquad \forall |x|<1.
\end{equation}
Then we have the following
\begin{theorem}\label{th-2d-1d}
Let be $\theta:=(\theta_1,\theta_2)\in ]0,1[^2$ and $f$ a locally continuous function on the square satisfying (\ref{hpf}) for the Jacobi weights (\ref{ui}).
For any $\n=(n_1,n_2)\in\NN^2$ and $\m=(m_1, m_2)\in \NN^2$ with {\color{black} $n_i\sim m_i$,} $i=1,2$, let $\V_\n^\m f$ be the bivariate VP filtered polynomial associated with the bivariate Chebyshev weight $v_1(x)v_2(y)$ interpolating $f$ at the Chebyshev grid $\X^\n=\{(x_k,y_h)\}_{k,h}$, being $x_k=x_k^{n_1}(v_1)$ ($k=1,\ldots, n_1$) and $y_h=x_h^{n_2}(v_2)$ ($h=1,\ldots, n_2$). If the couples of weights $w=v_i$ and $u=u_i$, for $i=1,2$ satisfy the bounds (\ref{hp-u}) then the error estimates
\begin{eqnarray}\label{eq-2d-1d}
\left|f(x,y)-\V_\n^\m f(x,y)\right|u_1(x)u_2(y)&\le& \C u_2(y) E_{n_1-m_1}(\overline{\overline{f}}_y)_{u_1}\\
\nonumber
&+& \C \sup_{1\le k\le n_1} u_1(x_k)E_{n_2-m_2}\left(\overline{f}_{x_k}\right)_{u_2}\\ [.15in]
\label{eq-2d-1d-bis}
\left|f(x,y)-\V_\n^\m f(x,y)\right|u_1(x)u_2(y)&\le& \C u_1(x) E_{n_2-m_2}(\overline{{f}}_x)_{u_2}\\
\nonumber
&+& \C \sup_{1\le h\le n_2} u_2(y_h)E_{n_1-m_1}\left(\overline{\overline{f}}_{y_h}\right)_{u_1}
\end{eqnarray}
hold for all $(x,y)\in [-1,1]^2$, with $C\ne \C(x,y,\n, f)$.
\end{theorem}
The proof of this theorem can be achieved similarly to \cite[Corollary 3.3]{paperoAMC} and for brevity it has been omitted.

We conclude, observing that from (\ref{eq-2d-1d}) several estimates can be deduced by considering the several classes of smoothness of $f(x,y)$ considered as a function of a single (the first or the second) variable. \newline
For instance, in the case that, for some $r,s>0$, $\overline{f}_x\in Z_r(u_2)$ and $\overline{\overline{f}}_x\in Z_s(u_1)$ hold true for all $|x|<1$, by means of (\ref{eq-2d-1d}) and (\ref{Favard-Zig}), for all $|x|,|y|\le 1$, we get
\begin{equation}\label{err-Zig-2d}
\left|f(x,y)-\V_\n^\m f(x,y)\right|u_1(x)u_2(y)\le \frac \C{n_1^s}+\frac \C{n_2^r},\qquad \C\ne\C(x,y,\n).
\end{equation}
\section{On simultaneous approximation }
In this section we discuss about the approximation of the $r$--th derivative ( $r\ge 1$) of a function $f$, namely $f^{(r)}$, in the case that it is completely unknown or not easily  computable and the only data at our disposal are the values of $f$ at Chebyshev nodes (\ref{zeros}).

In this case we can consider the approximation $f^{(r)}\approx ({V}_n^m f)^{(r)}$.  Indeed, by (\ref{VP}), we have
\begin{equation}\label{VP-r}
({V}_n^m f)^{(r)}(x)=\sum_{k=1}^n f(x_k^n)\left(\Phi_{n,k}^m\right)^{(r)}(x).
\end{equation}
In the following proposition an easy formula for the derivatives of the fundamental  VP polynomials is given.
\begin{proposition}\label{prop-der}
For any Chebyshev weight $w=v^{\alpha,\beta}$ and for all positive integers $n>m$ and $r\in\NN$, we have
\begin{equation}\label{fi-r}
\left(\Phi_{n,k}^m\right)^{(r)}(x)=\lambda_k^n\sum_{j=r}^{n-1} p_j(w,x_k^n)(q_{n,j}^m)^{(r)}(w,x),\qquad k=1,\ldots, n,
\end{equation}
{\color{black} and the derivatives $(q_{n,j}^m)^{(r)}(w,x)$ of the polynomial $q_{n,j}^m(w)$ at any $|x|\le 1$ are related to Jacobi orthonormal polynomials w.r.t. the Jacobi weight $(w\varphi^{2r})(x):=w(x)(1-x^2)^r$ as follows
\begin{equation}\label{qDER-basis}
(q_{n,j}^m)^{(r)}(w,x)=\left\{\begin{array}{ll}
\eta_{j}^r p_{j-r}(w\varphi^{2r},x)&\mbox{if}\quad 0\le j\le n-m,\\ [.1in]
\gamma_{n,j}^m\eta_{j}^r p_{j-r}(w\varphi^{2r},x)-\gamma_{n, 2n-j}^m\eta_{2n-j}^r p_{2n-j-r}(w\varphi^{2r},x)
&\mbox{if}\quad n\hspace{-.05cm}-\hspace{-.05cm}m<j<n,
\end{array}\right.
\end{equation}
where we set
\begin{equation}\label{eta-r}
\eta^{r}_n: =\eta^r_n(w)=\left(\prod_{k=0}^{r-1}(n-k)\times
 \prod_{k=1}^{r}(n+k+\chi_w)\right)^\frac 1 2,\quad \chi_w:=\al+\be,
\end{equation}
$\{\gamma_{n,j}^m\}_j$ are given by (\ref{gamma-q}) and we agree that $p_k(w\varphi^{2r},x)\equiv 0$ whenever $k<0$.}
\end{proposition}
{\it Proof of Proposition \ref{prop-der}}

Recalling that for any orthonormal Jacobi polynomial we have
\[
p_n'(v^{\alpha,\beta},x)=\sqrt{n(n+\alpha+\beta+1)}
p_{n-1}(v^{\alpha+1,\beta+1},x),\qquad n\ge 1,
\]
in the particular case  $|\al|=|\be|=\frac 1 2$,  we get
\begin{equation}\label{derivate_pn}
p_n^{(r)}(w,x)=\eta_{n}^r(w)p_{n-r}(w\varphi^{2r},x),\qquad n\ge r ,\qquad \forall r\in\NN,
\end{equation}
with $\eta_n^r$ given by (\ref{eta-r}).

Hence, the statement follows from (\ref{q-basis}) and (\ref{derivate_pn}).
\Proofend
\vspace{.1in}\newline
By formulas (\ref{VP-r}) and (\ref{fi-r}) we deduce that the computation of the derivatives of VP filtered interpolating polynomials can be performed by means of fast Fourier Discrete block transforms, similarly to the basic VP approximation. More precisely, for any $r\in\NN$ and any point set $\{z_\ell\}_{\ell=1,\ldots,N}$, setting
\[
\V^{(r)}:=\left[ \left(V_n^m f\right)^{(r)}(z_1),\ldots, \left(V_n^m f\right)^{(r)}(z_N) \right]
\]
and
\[
\mathcal{Q}^{(r)}:=\left[(q_{n,j}^m)^{(r)}(w,z_\ell)\right]_{0\le j< n,\ 1\le \ell \le N},
\]
we get
\begin{equation}\label{VPder-comput}
\V^{(r)} = \f \cdot\Lambda\cdot {\Phi}^{(r)},
\end{equation}
where $\f:=\left[f(x_1^n),\ldots, f(x_n^n)\right]$, $\Lambda:=\diag\left[\sqrt{  \lambda_k^n}\right]_{1\le k\le n}$ and we have
\begin{eqnarray}
\label{fi1-r} w=w_1 &\Longrightarrow&
{\Phi}^{(r)}= \IDC(\mathcal{Q}^{(r)}),\\
\label{fi2-r} w=w_2 &\Longrightarrow&
{\Phi}^{(r)}=\sqrt{\frac{n+1}{2}}\IDS(\mathcal{Q}^{(r)}),\\
\label{fi3-r} w=w_3 &\Longrightarrow&
{\Phi}^{(r)}=\frac{1}{\sqrt{2n+1}} \GDC(\mathcal{Q}^{(r)}),\\
\label{fi4-r} w=w_4 &\Longrightarrow&
{\Phi}^{(r)}=\frac{1}{\sqrt{2n+1}} \GDS(\mathcal{Q}^{(r)}),
\end{eqnarray}
being the previous transforms defined in Subsection 2.1.

In the sequel we are going to investigate the quality of the approximation  $f^{(r)}\approx (V_n^mf)^{(r)}$.

In the case that $f$ belongs to some Sobolev space, the following error estimate can be deduced from Theorem \ref{th-Sob}.
\begin{corollary}
Let the bounds (\ref{hp-u}) be satisfied and let $n,m\in\NN$ be such that {\color{black} $n\sim m$.} For any $r\in\NN$, and for any $f\in W_s(u)$  with $s\ge r$
\begin{equation}\label{err-simult}
\|(f-V_n^m f)^{(r)}u\varphi^r\|\le  \frac{\C}{n^{s-r}}\ \|f\|_{W_s(u)}, \qquad \C\ne \C(n,f).
\end{equation}
\end{corollary}
Finally, in the case of the 1st kind Chebyshev weight, we prove the following  unweighted  error estimate
\begin{theorem}\label{simultanea}
Let be $w=w_1$. {\color{black} For any $n,m\in\NN$ such that $n\sim m$ and for each }$f\in C^{q}([-1,1])$, the associated VP filtered interpolating polynomial $V_n^mf$ and its derivatives of order $r\le q/2$, satisfy the following uniform error estimate
\begin{equation}\label{err-puntuale}
\|f^{(r)}-\left(V_n^m f \right)^{(r)}\|\le \frac{\C}{n^{q-2r}}E_{n-m-q}(f^{(q)}), \qquad \C\ne \C(n,f).
\end{equation}
%for any $0\le 2r\le q.$
\end{theorem}

\vspace{0.3cm}

{\it Proof of Theorem \ref{simultanea}}. \
By the Gopengauz theorem \cite{gopengauz}, for any $f\in C^{q}([-1,1])$ there exists a polynomial $P\in \PP_{n-m}$ satisfying, for any $0\le r\le q$, the following estimate
\begin{equation}\label{gop}|f^{(r)}(x)-P^{(r)}(x)|\le \C   \left(\frac{\sqrt{1-x^2}}{n-m} \right)^{q-r}E_{n-m-q}(f^{(q)}), \ \ \forall x\in [-1,1]. \end{equation}
By using this polynomial and recalling that $V_n^m P =P$, we have
\begin{eqnarray}
\nonumber
\left|\left(f-V_n^m f \right)^{(r)}(x)\right|&=&
\left|\left(f-P \right)^{(r)}(x)+ \left(V_n^mf -V_n^mP\right)^{(r)}(x)\right|
\\
\nonumber&\le& \left|f^{(r)}(x)-P^{(r)}(x)\right|+
\left|\left(V_n^m (f-P) \right)^{(r)}(x)\right|\\
&=:&I_1(x)+I_2(x).\label{I1+I2}
\end{eqnarray}
By (\ref{gop}) and by taking into account $(n-m)=(1-\theta) n$, we have
\begin{equation}\label{I1}
I_1(x)\le \C   \frac{E_{n-m-q}(f^{(q)})}{n^{q-r}}, \qquad\forall |x|\le 1, \qquad \C\ne \C(n,x,f).\end{equation}
In order to estimate $I_2(x)$, we recall the following Markov inequality %
\cite[p.233, Th. 5.1.8]{borwein}
\begin{equation}\label{Markov}
\|Q^{(r)}\|\le \C \nu^{2r}\|Q\|,\quad \forall Q\in\PP_\nu, \qquad \C\ne\C(Q, \nu).
\end{equation}
Moreover, since the couple of weights $w=w_1$ and $u=1$ (i.e. $\gamma=0=\delta$) satisfy the bounds (\ref{hp-u}), by Theorem \ref{th-C0} we have that $V_n^m: C^0\rightarrow C^0$ is an uniformly bounded map, i.e.
\begin{equation}\label{bound-VP}
\|V_n^m g\|\le \C \|g\|,\quad \forall g\in C[-1,1], \qquad \C\ne\C(n,g).
\end{equation}
Hence, by (\ref{Markov}), (\ref{bound-VP}) and (\ref{gop}), we get
\begin{eqnarray}\nonumber
I_2(x)&:=& \left|\left(V_n^m (f-P) \right)^{(r)}(x)\right|\le \C  (n+m-1)^{2r}\|V_n^m (f-P)\|\\
\nonumber
&\le&\C n^{2r} \|V_n^m (f-P)\|\le \C n^{2r}\|f-P\|\\
&\le& \C \frac{n^{2r}}{n^q}E_{n-m-q}(f^{(q)}), \qquad\forall |x|\le 1, \qquad \C\ne \C(n,x,f).
\label{I2}
\end{eqnarray}
Therefore the conclusion follows by combining  (\ref{I1}) and (\ref{I2}) with (\ref{I1+I2}).
\Proofend

\section{On the choice of the parameter $m$}
One of the main advantages of VP filtered interpolation consists of improving the approximation by keeping fixed the data (i.e. the number $n$ of interpolation nodes) and suitably modulating the action ray $m$ of the VP filter. In several papers  \cite{woula_NA, ThBa, lncs_nostro, paperoAMC, filtered1D} this feature has been noticed in connection with the uniform boundedness of the Lebesgue constants that are defined as the operator norm of the map $V_n^m:C^0_u\rightarrow C^0_u$, i.e.
\[
\|V_n^m\|_{C^0_u\rightarrow C^0_u}:= \sup_{f\in C^0_u}\frac{\|(V_n^mf)u\|}{\|fu\|}, \qquad 0<m<n.
\]
By Theorem \ref{th-C0}, for suitable choices of the weight $u$, these constants are uniformly bounded w.r.t. $n$, if we take $m=\theta n$ for an arbitrarily fixed $\theta\in ]0,1[$.
In this section, we are going to provide some additional hints for an optimal choice of the degree--parameter $m$, once fixed the number $n$ of data.

We limit our analysis to the case that $u:=\sqrt{w\varphi}$, where $w$ is the Chebyshev weight associated with $V_n^mf$. In this way, the bounds (\ref{hp-u}) required by Theorem \ref{th-C0} are certainly satisfied for any $w$.

Note that if we apply (\ref{inva}) to the polynomial $P^*\in \PP_{n-m}$ s.t. $\|(f-P^*)u\|=E_{n-m}(f)_u$, then we get
\begin{eqnarray*}
\|(V_n^mf-f)u\|&=& \|(V_n^m (f-P^*)+P^*-f)u\|\\
&\le& \left[\|V_n^m\|_{C^0_u\rightarrow C^0_u}+1\right]\|(f-P^*)u\|\\
&=& \left[\|V_n^m\|_{C^0_u\rightarrow C^0_u}+1\right]\ E_{n-m}(f)_u.
\end{eqnarray*}
Hence, a possible value for the constant $\C>0$ in (\ref{err-VP-C0}) is any upper bound for the constant
\begin{equation}\label{Cnm}
c(n,m):=1+\|V_n^m\|_{C^0_u\rightarrow C^0_u},
\end{equation}
which is obviously independent of $f$ and, according with Theorem \ref{th-C0}, it is also independent of $n$ and $m$ whenever we take $n,m\in\NN$ s.t. {\color{black} $n\sim m$.}

Recently, in \cite[Table 1]{filtered1D} for different choices of the weight functions $u$ and $w$,  the authors have numerically computed
\[
L_\theta:=\sup_{\footnotesize \color{black}
\begin{array}{c}n,m\in\NN \\ m=\theta n
\end{array}}\|V_n^m\|_{C^0_u\rightarrow C^0_u}, \qquad\quad \mbox{for $\theta=0.1,\ 0.2,\ldots, 0.9$},
\]
and  all the performed numerical tests displayed a decreasing behavior of $L_\theta$ w.r.t. $\theta$.

Such a behavior is supported by the following theorem.
\begin{theorem}\label{th-theta}
For any pair of positive integers $n>m$, let $V_n^m:f\rightarrow V_n^mf$ be the VP operator corresponding to a Chebyshev weight $w$. For $u=\sqrt{w\varphi}$, we have
\begin{equation}\label{LC-theta}
\|V_n^m\|_{C^0_u\rightarrow C^0_u}\le C_w\sqrt{\frac {n_w}{m}}\left[1+\frac{2\pi m_w}{n_w}\right],
\end{equation}
where $C_w=1$ unless $w=w_1$ when it is $C_{w_1}=2$,
\begin{equation}\label{nw}
n_w:=\left\{\begin{array}{ll}
n & \mbox{if $w=w_1$,}\\ [.1in]
n+1 & \mbox{if $w=w_2$,}
\\ [.1in]
\displaystyle
\frac{2n+1}{2} & \mbox{if $w\in \{w_3,w_4\},$}
\end{array}\right.
\end{equation}
and
\[
m_w:= \left\{\begin{array}{ll}
n+m & \mbox{if $w\in\{w_1, w_2\}$,}\\ [.1in]
2(n+m)-1 & \mbox{if $w\in\{w_3, w_4\}$}.
\end{array}\right.
\]
\end{theorem}
Before giving the proof of this theorem let us make some comments.

{\color{black} By virtue of (\ref{LC-theta}) it turns out that, according with Theorem \ref{th-C0}, the constant $c(n,m)$ in (\ref{Cnm}) is uniformly bounded w.r.t. $n$ whenever we take the parameter $m$ such that $m<n\le \C m$ with $\C\ne \C(n,m)$.
A usual choice is to fix a rational parameter $\theta\in ]0,1[$ and take the integers $m=\theta n$. In this particular case, the observed decreasing behavior of the Lebesgue constants $L_\theta$ is supported by the right hand side in (\ref{LC-theta}), which is a decreasing function of $m<n$.}

Finally, we remark that even if any upper bound of the Lebesgue constant gives an upper bound of the constant $c(n,m)$,  we have to take into account that this is only one factor of the error estimate
\[
\|(V_n^mf-f)u\|\le c(n,m) E_{n-m}(f)_u, \qquad c(n,m):=1+\|V_n^m\|_{C^0_u\rightarrow C^0_u}.
\]
The other factor, $E_{n-m}(f)_u$, depends on the smoothness of $f$, but it depends on $n$ and $m$ too. For instance, supposed that $f\in Z_s(u)$ for some $s>0$, then we have
\begin{equation}\label{eq-s}
E_{n-m}(f)_u\le \frac{\|f\|_{Z_s(u)}}{(n-m)^s},
\end{equation}
which is an increasing function of $m$.

Hence, supposed to know the degree $s>0$ of smoothness of $f$, for an optimal choice of the parameter $m$ with a certain fixed $n$,  we have to balance between the decreasing function of $m$ at the right hand side of (\ref{LC-theta}) and the increasing function of $m$ at the right hand side of (\ref{eq-s}), since we have to take into account that
\begin{equation}\label{VP-err-nm}
\|(V_n^mf-f)u\|\le \left(1+\C_w\left[\sqrt{\frac{n_w}m}+\frac{2\pi m_w}{\sqrt{m\ n_w}}\right]\right)\frac{\|f\|_{Z_s(u)}}{(n-m)^s}, \qquad \forall f\in Z_s(u).
\end{equation}
In conclusion, let us prove Theorem \ref{th-theta}. To this aim, we need to state a preliminary lemma that can be useful in other contexts too.
\begin{lemma}\label{lemmarci}
For any $n,\nu\in\NN$, let $Q$ be an even or odd trigonometric polynomial of degree at most $\nu$ and let $t_k^n=t_k^n(w)$, $k=1,\ldots,n$, be given by (\ref{zeros}). We have
\begin{equation}\label{Marci}
\sum_{k=1}^n |Q(t_k^n)|\Delta t_k^n\le\left(1+\frac{2\pi\nu}{n_w}\right) \int_0^\pi |Q(\tau)|d\tau,
\end{equation}
where $\Delta t_k^n:= (t_k^n-t_{k-1}^n)$ for $k=1,\ldots,n$ ( with $t_{0}^n:=0$) and $n_w$ has been  defined in (\ref{nw}).
\end{lemma}
{\it Proof of Lemma \ref{lemmarci}}

By the Mean Value Theorem, we have
\begin{equation}\label{meanVT}
\int_{t_{k-1}^n}^{t_k^n}|Q(\tau)|d\tau =|Q(\xi_k)|\Delta t_k^n, \qquad t_{k-1}^n< \xi_k< t_{k}^n,\qquad k=1,\ldots, n.
\end{equation}
Using the points $\{\xi_k\}_{k=1,..,n}$ and  (\ref{meanVT}), we note that
\begin{eqnarray*}
\sum_{k=1}^n |Q(t_k^n)|\Delta t_k^n&\le&
\sum_{k=1}^n |Q(\xi_k)|\Delta t_k^n+\sum_{k=1}^n \left|Q(t_k^n)-Q(\xi_k)\right|\Delta t_k^n\\
&=& \sum_{k=1}^n \int_{t_{k-1}^n}^{t_{k}^n}|Q(\tau)|d\tau +
\sum_{k=1}^n \left|\int^{t_k^n}_{\xi_k} Q^\prime (\tau)d\tau\right|\Delta t_k^n\\
&\le& \int_0^{\pi} |Q(\tau)|d\tau +
\left(\max_{1\le k\le n}\Delta t_k^n\right)\sum_{k=1}^n\int_{t_{k-1}^n}^{t_{k}^n} \left|Q^\prime (\tau)\right|d\tau\\
&\le& \int_{0}^\pi |Q(\tau)|d\tau +
\left(\max_{1\le k\le n}\Delta t_k^n\right)\int_{0}^{2\pi} \left|Q^\prime (\tau)\right|d\tau.
\end{eqnarray*}
On the other hand, we observe that
\begin{equation}\label{Dtk}
\Delta t_k^n=\left\{\begin{array}{ll}
\displaystyle \frac \pi {n} \quad \left(\frac \pi{2n} \ \mbox{if $k=1$}\right) & \mbox{if $w=w_1$}\\ [.1in]
\displaystyle\frac \pi{n+1} & \mbox{if $w=w_2$}
\\ [.1in]
\displaystyle\frac{2 \pi}{2n+1} & \mbox{if $w=w_3$ or $w=w_4$}
\end{array}\right.\qquad k=1,\ldots,n,
\end{equation}
which yields
\begin{equation}\label{Dtk-max}
\max_{1\le k\le n}\Delta t_k^n= \frac\pi{n_w}.
\end{equation}
 Hence, by applying the Bernstein inequality \cite[Th. 5.1.4]{borwein}
\begin{equation}\label{Ber}
\int_{0}^{2\pi} \left|Q^\prime (\tau)\right|d\tau\le \nu \int_{0}^{2\pi} \left|Q (\tau)\right|d\tau,
\qquad \nu=\deg (Q),
\end{equation}
and recalling that $Q$ is an even or odd trigonometric polynomial, we obtain the statement as follows
\begin{eqnarray*}
\sum_{k=1}^n |Q(t_k^n)|\Delta t_k^n&\le&
 \int_{0}^\pi |Q(\tau)|d\tau +
\frac\pi{n_w}\int_{0}^{2\pi} \left|Q^\prime (\tau)\right|d\tau
\\
&\le&
 \int_{0}^{\pi} |Q(\tau)|d\tau +
\frac{\pi\nu}{n_w}\int_{0}^{2\pi} \left|Q (\tau)\right|d\tau
\\
&=& \left[1+ \frac{2\pi\nu}{n_w}\right] \int_{0}^\pi |Q(\tau)|d\tau .
 \end{eqnarray*}\Proofend

{\it Proof of Theorem \ref{th-theta}}

By (\ref{VP}) we easily derive that
\begin{equation}\label{eq-LC1}
\|V_n^m\|_{C^0_u\rightarrow C^=_u}=\sup_{|x|\le 1}\sum_{k=1}^n \left|\Phi_{n,k}^m(x)\right|\frac{u(x)}{u(x_k^n)}.
\end{equation}
On the other hand, in proving \cite[Proposition 2.1]{filtered1D} the following trigonometric form of the fundamental VP polynomials has been deduced by means of classical trigonometric identities (see  \cite[pp. 13--14]{filtered1D})

$\bullet${\it Case $w=w_1$}
\begin{eqnarray*}
\Phi_{n,k}^m(\cos t)&=&
\frac 1{4 nm}\left[\frac{\sin[n(t-t_k)]\sin[m(t-t_k)]}{\sin^2[(t-t_k)/2]}\right.\\
&&\left. + \frac{\sin[n(t+t_k)]\sin[m(t+t_k)]}{\sin^2[(t+t_k)/2]}\right].
\end{eqnarray*}
$\bullet${\it Case $w=w_2$}
\begin{eqnarray*}
\Phi_{n,k}^m(\cos t)&=&
\frac{\sin t_k}{4m(n+1)\sin t}
\left[\frac{\sin[(n+1)(t-t_k)]\sin[m(t-t_k)]}{\sin^2[(t-t_k)/2]}\right.\\
&&\left.- \frac{\sin[(n+1)(t+t_k)]\sin[m(t+t_k)]}{\sin^2[(t+t_k)/2]}\right].
\end{eqnarray*}
$\bullet$ {\it Case $w=w_3$}
\begin{eqnarray*}\Phi_{n,k}^m(\cos t)&=&
\frac{\cos [t_k/2]}{2m(2n+1)\cos [t/2]}
\left[\frac{\sin[(2n+1)(t-t_k)/2]\sin[m(t-t_k)]}{\sin^2[(t-t_k)/2]} +\right.\\ [.1in]
&& \left.+ \frac{\sin[(2n+1)(t+t_k)/2]\sin[m(t+t_k)]}{\sin^2[(t+t_k)/2]}\right].
\end{eqnarray*}
$\bullet${\it Case $w=w_4$}
\begin{eqnarray*}
\Phi_{n,k}^m(\cos t)&=&
\frac{\sin [t_k/2]}{2m(2n+1)\sin [t/2]}
\left[\frac{\sin[(2n+1)(t-t_k)/2]\sin[m(t-t_k)]}{\sin^2[(t-t_k)/2]} +\right.\\ [.1in]
&& \left.- \frac{\sin[(2n+1)(t+t_k)/2]\sin[m(t+t_k)]}{\sin^2[(t+t_k)/2]}\right].
\end{eqnarray*}
Hence, setting
\[
\Psi(N,M, \tau):= \frac{\sin [N\tau/2]\sin[M\tau/2]}{\sin^2[\tau/2]}
\]
and recalling that $u=\sqrt{w\varphi}$, for any $x=\cos t$, $t\in [0,\pi]$, we get
\[
\sum_{k=1}^n \left|\Phi_{n,k}^m(x)\right|\frac{u(x)}{u(x_k^n)}=\sum_{k=1}^n \left\{
\begin{array}{ll}
\displaystyle \frac{\left|\Psi\left(2n,\ 2m,\ t-t_k^n\right)+\Psi\left(2n,\ 2m,\ t+t_k^n\right)\right|}{4mn} & \mbox{if $w=w_1$},\\ [.15in]
\displaystyle \frac{\left|\Psi\left(2(n+1),\ 2m,\ t-t_k^n\right)-\Psi\left(2(n+1),\ 2m,\ t+t_k^n\right)\right|}{4m(n+1)} & \mbox{if $w=w_2$},\\ [.15in]
\displaystyle\frac{\left|\Psi\left(2n+1,\ 2m,\ t-t_k^n\right)+\Psi\left(2n+1,\ 2m,\ t+t_k^n\right)\right|}{2m(2n+1)} & \mbox{if $w=w_3$},\\[.15in]
\displaystyle\frac{\left|\Psi\left(2n+1,\ 2m,\ t-t_k^n\right)-\Psi\left(2n+1,\ 2m,\ t+t_k^n\right)\right|}{2m(2n+1)} & \mbox{if $w=w_4$}.
\end{array}\right.\]
In order to write this formula in a more convenient form, we set
\begin{equation}\label{Q-def}
Q(\tau):=\left\{
\begin{array}{ll}
\Psi\left(2n,\ 2m,\ t-\tau\right)+\Psi\left(2n,\ 2m,\ t+\tau\right)& \mbox{if $w=w_1$},\\ [.15in]
\Psi\left(2(n+1),\ 2m,\ t-\tau\right)-\Psi\left(2(n+1),\ 2m,\ t+\tau\right) & \mbox{if $w=w_2$},\\ [.15in]
\Psi\left(2n+1,\ 2m,\ t-\tau\right)+\Psi\left(2n+1,\ 2m,\ t+\tau\right)& \mbox{if $w=w_3$},\\[.15in]
\Psi\left(2n+1,\ 2m,\ t-\tau\right)-\Psi\left(2n+1,\ 2m,\ t+\tau\right) & \mbox{if $w=w_4$},
\end{array}\right.\end{equation}
so that, taking into account (\ref{Dtk}), the previous formula becomes
\[
\sum_{k=1}^n \left|\Phi_{n,k}^m(x)\right|\frac{u(x)}{u(x_k^n)}
\left\{\begin{array}{ll}
\displaystyle
=\frac 1{4m\pi}\sum_{k=1}^n |Q(t_k^n)|\Delta t_k^n & \mbox{if $w\ne w_1$}\\
\displaystyle
\le \frac 1{2m\pi}\sum_{k=1}^n |Q(t_k^n)|\Delta t_k^n & \mbox{if $w=w_1$}
\end{array}\right.
\]
i.e., for any Chebyshev weight $w$, we have
\[
\sum_{k=1}^n \left|\Phi_{n,k}^m(x)\right|\frac{u(x)}{u(x_k^n)}
\le \frac {C_w}{4m\pi}\sum_{k=1}^n |Q(t_k^n)|\Delta t_k^n.
\]
Now, in the case $w\in\{w_1, w_2\}$, we observe that $Q(\tau)$ is an even (if $w=w_1$) or odd (if $w=w_2$) trigonometric polynomial of degree at most $\nu=n+m$. Hence, by  Lemma \ref{lemmarci} we get
\[
\sum_{k=1}^n \left|\Phi_{n,k}^m(x)\right|\frac{u(x)}{u(x_k^n)}\le \frac {C_w}{4m\pi}
\sum_{k=1}^n |Q(t_k^n)| \Delta t_k^n\le \frac {C_w}{4m\pi}\left[ 1+\frac{2\pi \nu}{n_w}\right]
\int_0^\pi |Q(\tau)|d\tau
\]
i.e., since $|Q(-\tau)|=|Q(\tau)|$, we have
\begin{equation}\label{int-Q}
\sum_{k=1}^n \left|\Phi_{n,k}^m(x)\right|\frac{u(x)}{u(x_k^n)}\le
\frac {C_w}{8m\pi}\left[1+\frac{2\pi \nu}{n_w}\right]
\int_{-\pi}^{\pi} |Q(\tau)|d\tau, \qquad \nu=n+m.
\end{equation}
In the case $w\in \{w_3,w_4\}$, $Q(\tau)$ is not a trigonometric polynomial. Nevertheless, by following the same steps of the proof of Lemma \ref{lemmarci}, from the Mean Value Theorem and (\ref{Dtk-max})  we derive that
\[
\sum_{k=1}^n |Q(t_k^n)| \Delta t_k^n\le  \int_{0}^\pi |Q(\tau)|d\tau +
\frac\pi{n_w}\int_{0}^{\pi} \left|Q^\prime (\tau)\right|d\tau.
\]
On the other hand, $C:= Q(2\tau)$ is an even (if $w=w_3$) or odd (if $w=w_4$) trigonometric polynomial of degree $2\nu -1$, $\nu=n+m$, w.r.t. the variable $\tau$. Hence, by the variables change $\tau=2s$ and the Bernstein inequality (\ref{Ber}) applied to $\tilde Q$, we get
\[
\int_0^{\pi}|\left|Q^\prime (\tau)\right|d\tau\le \int_0^{2\pi}\left|\tilde Q^\prime (s)\right|ds\le
(2\nu-1) \int_0^{2\pi}\left|\tilde Q (s)\right|ds= 2(2\nu-1) \int_0^{\pi}|Q(\tau)|d\tau.
\]
Summing up, if $w\in\{w_3,w_4\}$ then we have
\[
\sum_{k=1}^n |Q(t_k^n)| \Delta t_k^n\le \left[1+\frac{2\pi (2\nu-1)}{n_w}\right] \int_{0}^\pi |Q(\tau)|d\tau,
\]
and consequently
\begin{equation}\label{int-Q-1}
\sum_{k=1}^n \left|\Phi_{n,k}^m(x)\right|\frac{u(x)}{u(x_k^n)}\le
\frac {C_w}{8m\pi}\left[1+\frac{2\pi (2\nu-1)}{n_w}\right]
\int_{-\pi}^{\pi} |Q(\tau|d\tau, \qquad \nu:=n+m.
\end{equation}
Finally, we observe that for any $t\in [0,\pi]$ we have
\begin{eqnarray*}
\int_{-\pi}^\pi \left|\psi(N,M,t\pm \tau)\right|d\tau &=&
\int_{-\pi}^\pi \left|\psi(N,M,\tau)\right|d\tau=
\int_{-\pi}^\pi \left|\frac{\sin[N\tau/2]\sin[M\tau/2]}{\sin^2[\tau/2]}\right|d\tau\\
&\le& \left(\int_{-\pi}^\pi \left|\frac{\sin^2[N\tau/2]}{\sin^2[\tau/2]}\right|d\tau\right)^\frac 12
\left(\int_{-\pi}^\pi \left|\frac{\sin^2[M\tau/2]}{\sin^2[\tau/2]}\right|d\tau\right)^\frac 12,
\end{eqnarray*}
and taking into account that (see e.g. \cite[p.26]{mastromilobook})
\[
\frac 1{2\pi}\int_{-\pi}^{\pi}\frac{\sin^2[K\tau/2]}{\sin^2[\tau/2]}d\tau =K,\qquad \forall K\in\NN,
\]
we get
\begin{equation}\label{int-psi}
\int_{-\pi}^\pi \left|\psi(N,M,t\pm \tau)\right|d\tau\le 2\pi\sqrt{M\cdot N}, \qquad \forall N,M\in\NN.
\end{equation}
Consequently, by means of (\ref{Q-def}) and (\ref{int-psi}), for any $w$, we have
\[
\frac {1}{8m\pi}\int_{-\pi}^\pi|Q(\tau)|d\tau \le \frac 1{4m\pi}\int_{-\pi}^\pi |\Psi(2n_w,2m,\tau)|d\tau
\le  \frac{\sqrt{m\cdot n_w }}{m}= \sqrt{\frac{ n_w} m},
\]
which, combined with the estimates (\ref{int-Q}) and (\ref{int-Q-1}), yields the statement.
\Proofend

\noindent\textbf{Acknowledgments.}
This research has been accomplished within Rete ITaliana di Approssimazione (RITA) and partially supported by the GNCS-INdAM funds 2019, project ``Discretizzazione di misure, approssimazione di operatori integrali ed applicazioni''.

%\bibliographystyle{plain}
%
%\bibliography{biblio}

{\small{\it{
Donatella Occorsio, {Department of Mathematics, Computer Science and Economics,\\ University of Basilicata,\\  Via dell'Ateneo Lucano 10, 85100 Potenza, Italy }. \\ donatella.occorsio@unibas.it.

{{Woula Themistoclakis, C.N.R. National Research Council of Italy, IAC Institute for Applied Computing ``Mauro Picone'',\\ Via P. Castellino, 111, 80131 Napoli, Italy. \\ {woula.themistoclakis@cnr.it.}}}
  }}}

\end{document}